\documentclass[12pt]{article}
\usepackage{amssymb,amsmath}
\newtheorem{theorem}{Theorem}[section]
\newtheorem{definition}{Definition}[section]

\newtheorem{remark}{Remark}[section]

\begin{document}
\begin{center}
{\large \bf Fixed points of left reversible semigroup of isometry mappings in  Banach spaces}\\
\vspace{0.5cm}

{ S. Rajesh}

{Department of Mathematics}

{ Indian Institute of Technology Tirupati, India - 517506 }

{e-mail: srajeshiitmdt@gmail.com}
\end{center}

\begin{abstract}
In this paper, we prove the existence of a common fixed point in $C(K)$, the  Chebyshev center of K, for a left reversible semigroup of isometry mappings. This existence result improves the results obtained by Lim et al. and Brodskii and Milman.  
\end{abstract} 

 \section{Introduction} 
Let $K$ be a nonempty weakly compact convex set in a Banach space $X$. A map $T:K \rightarrow X$ is said to be isometry if $d(Tx, Ty) = d(x, y)$ for all $x, y \in K$. The notion of 
normal structure, which is introduced in \cite{Brod_1948}, is defined as follows:
\begin{definition}\label{P:D:1}\cite{Brod_1948,GK_Book}
A nonempty bounded convex set $K$ in a Banach space $X$
is said to have normal structure if for every nonempty convex
set $C \subseteq K$ with $\delta(C) > 0$ has a point $x \in C$
such that $r(x, C) < \delta(C)$, where $r(x, C) = \sup\{\|x - y\|: y \in C\}$
and $\delta(C) = \sup\{\|z - y\| : z, y \in C\}$.

 Define $r(K) = \inf\{r(x, K): x \in K\}$ and $C(K) = \{x \in K: r(x, K) = r(K)\}$.
 Then the set $C(K)$ and the number $r(K)$ are called, respectively,
 the set of Chebyshev center of $K$ and the Chebyshev radius of $K$.
\end{definition}

The notion of UKK--norm is defined as follows:

\begin{definition} \cite{Bena_Book,GK_Book}
A Banach space $X$ is said to have uniformly Kadec-Klee (UKK) norm
if and only if for any $\epsilon > 0$, there exists $\delta > 0$ such that
\begin{center}
   $\{x_n\} \subseteq B[0, 1]$, $x_n$ converges weakly to $x_0$,
  and ${\it sep}\{x_n\}:=\inf\{\|x_n - x_m\|: n \neq m\} > \epsilon$,
\end{center}
imply that
\begin{center}
  $\|x_0\| \leq 1 - \delta$.
\end{center}
\end{definition}


\begin{definition}\cite{Cliff_1961,Mitc_1970}
A semigroup $S$ is called left reversible if for every pair of elements $a,b \in S$, there exists a pair $c, d \in S$ such that $ac=bd$. 
\end{definition}

Now, we state the basic common fixed point theorem for left reversible semigroup of nonexpansive maps.

\begin{theorem}\label{P:T:CFT_1}\cite{Mitc_1970}
Let $K$ be a nonempty compact convex set in a Banach space and $\mathfrak{F}$ 
be a left reversible semigroup of nonexpansive self-maps on $K$. Then $K$ contains a common fixed point of $\mathfrak{F}$.
\end{theorem}

\begin{theorem}\label{P:T:CFT_2}\cite{Lim_1974_1}
Let $K$ be a nonempty weakly compact convex set in a Banach space $X$ and assume that $K$ has normal structure. Let $S$ be a left reversible topological semigroup of nonexpansive, separately continuous actions on $K$. Then $K$ contains a common fixed point for $S$.
\end{theorem}

For the recent advancements in existence of common fixed points for various semigroups, one can refer to \cite{Lau_2002,Lau_2008,Lau_2010} and the reference therein.

The study of existence of fixed points in $C(K)$ is initiated by Brodskii and Milman \cite{Brod_1948}. In fact, Brodskii and Milman
proved: 

\begin{theorem}\label{P:T:Brod}\cite{Brod_1948}
  Let $K$ be a nonempty weakly compact convex subset of a Banach space $X$ and
  $\mathfrak{F} =\{T: K \rightarrow K: T$ is a surjective isometry mapping$\}$.
    Assume that $K$ has normal structure,
then there exists $x \in C(K)$ 
  such that $Tx = x$ for every $T \in \mathfrak{F}$.
\end{theorem} 

Motivated by this result (Theorem \ref{P:T:Brod}) of Brodskii and Milman and the fact
$T(C(K)) = C(K)$ whenever $T$ is a surjective isometry on $K$,
Lim et al. raised the following questions in \cite{Lim_2003}:

\noindent{\bf Question 1.} Let $T$ be an isometry on $K$ which is not surjective.
Does one still have $T(C(K)) \subseteq C(K)$ ?

\noindent{\bf Question 2.} Let $K$ be a weakly compact convex subset in a
Banach space $X$ and assume that $K$ has normal structure. Does there exist
a point in $C(K)$ which is fixed by every isometry from $K$ into $K$ ?

In case of uniformly convex Banach spaces, Lim et al. \cite{Lim_2003} affirmatively answered the above questions. 
Also, Lim et al. \cite{Lim_2003} established the next result:

\begin{theorem}\label{P:T:Lim_1}\cite{Lim_2003}
  Let $K$ be a nonempty weakly compact convex set in a 
   Banach space $X$ and $T$ be an isometry from $K$ into $K$. Assume that $K$ has normal structure. Then $T$ has a fixed point in $C(K)$.
\end{theorem}

In the setting of strictly convex Banach spaces, the authors in \cite{Raj_2015} proved: 
\begin{theorem}\cite{Raj_2015}
Let $K$ be a nonempty weakly compact convex set $K$ having normal structure in a strictly convex Banach space $X$ and $\mathfrak{F}$ be a commuting family of isometry self-mappings on $K$. Then $\mathfrak{F}$ has a common fixed point in the set of Chebyshev center, $C(K)$, of $K$.
\end{theorem}  

Moreover, the authors in \cite{Raj_2016_jnca} showed that the set of all Chebyshev center $C(K)$ need not be invariant under isometry mappings. 

In this paper, we prove that if $K$ is a nonempty weakly compact convex set in a Banach space $X$ and $\mathfrak{F}$ is a left reversible semigroup of isometry mappings on $K$ such that every sub-semigroup of $\mathfrak{F}$ is also left reversible, then there exists a point $x_0$ in $C(K)$, the set of Chebyshev center of $K$, such that $Tx_0 = x_0$ for all $T \in \mathfrak{F}$, whenever either $X$ has UKK--norm or $X$ is strictly convex and $K$ has normal structure.

\section{Common fixed point theorems}
Now, we state an interesting fact about isometry mappings on compact sets in metric spaces. 
\begin{theorem}\label{P:T:.1}\cite{Kuma_Book}
Let $K$ be a nonempty compact set in a metric space and $T: K \rightarrow K$ be an isometry map. Then $T(K) = K$.
\end{theorem}

\noindent{\bf Proof} Suppose $T(K) \subsetneqq K$. Then there exists $y \in K$ such that 
$T^{n}(y) \notin T^{n+1}(K)$ for $n = 0, 1, 2, ...$, where $T^0(y) = y.$ 
Since $T(K)$ is compact and $y \notin T(K)$, there exists $z \in T(K)$ such that $d(y, z) = dist(y, T(K))$, where $dist(y, T(K)) = \inf\{d(y, x): x \in T(K)\}$. 

Also, it is easy to see that $dist(T^n(y), T^{n+1}(K)) = dist(y, T(K))$, for all $n \in \mathbb{N}$. Hence $dist(T^n(y), T^{n+1}(K)) = d(y, z), $ for all $n \in \mathbb{N}$.

Since $\{T^n(y)\}$ is a sequence in the compact set $K$, $\{T^n(y)\}$ has a convergent subsequence, say $\{T^{n_k}(y)\}$. Suppose that $\{T^{n_k}(y)\}$ converges to $y_0 \in K$. Note that $T^{n}(K) \subseteq T^{n-1}(K)$, for $n \in \mathbb{N}$ and $\{T^{n+1}(y)\} \subseteq T^n(K)$ for all $n \in \mathbb{N}$. Hence $y_0 \in \bigcap_{n \in \mathbb{N}}T^n(K).$ 

As $y_0 \in T^n(K)$ for $n \in \mathbb{N}$, $d(T^n(y), y_0) \geq dist(T^n(y), T^{n+1}(K)) = d(y, z)$. Hence $\lim_{k \rightarrow \infty} d(T^{n_k}(y), y_0) \geq d(y, z) > 0$. But the sequence $\{d(T^{n_k}(y), y_0)\}$ converges to $d(y_0, y_0) = 0.$
This contradiction shows that $T(K) = K$.
\hfill{$\Box$}

\begin{theorem}\label{P:T:.2}
Let $K$ be a nonempty compact convex set in a Banach space $X$ and $\mathfrak{F}$ be a family of isometry self-maps on $K$. Then there exists $x_0 \in C(K)$ such that $Tx_0 = x_0$ for all $T \in \mathfrak{F}$. 
\end{theorem}

\noindent{\bf Proof}
Note that $K$ has normal structure, as $K$ is a compact convex set. Also, it follows from Theorem \ref{P:T:.1} that $T(K) = K$ for all $T \in \mathfrak{F}$. Hence, by Theorem \ref{P:T:Brod}, there exists $x_0 \in C(K)$ such that 
$Tx_0 = x_0$, for all $T \in \mathfrak{F}$.
\hfill{$\Box$}

\begin{theorem}\label{P:T:1} Let $K$ be a nonempty weakly compact convex set in a Banach space $X$ with a UKK--norm. Let $\mathfrak{F}$ be a left reversible topological semigroup of isometry, separately continuous actions on $K$. Assume that every sub-semigroup of $\mathfrak{F}$ is also left reversible. Then the Chebyshev center of $K$, $C(K)$, contains a common fixed point of $\mathfrak{F}$. 
\end{theorem} 

\noindent{\bf Proof}
It follows from Theorem \ref{P:T:Lim_1} that 
every map in $\mathfrak{F}$ has a fixed point in $C(K)$. Let $F_T = \{x \in C(K): Tx = x\}$ 
for $T \in \mathfrak{F}$. It is known \cite{Bena_Book,GK_Book} that if $X$ has a UKK--norm, then $C(K)$ is a nonempty compact subset of $K$. Hence it is enough to prove that the family $\mathfrak{S} = \{F_T: T \in \mathfrak{F}\}$ has the finite intersection property. 

Now, consider a finite subfamily $\{T_1, T_2, ... T_m\}$ of $\mathfrak{F}$. Take $m=3$, the same proof 
will go through for any $m \geq 4.$ Let $n \geq 3$ and $\sigma \in \mathfrak{S}_n$, where $\mathfrak{S}_n$ is the permutation group.  

Set $T_{\sigma(j)} = T_l$ if $\sigma(j) \equiv l (mod~3)$ and 
$T^{\sigma} = T_{\sigma(1)}\circ T_{\sigma(2)}\circ...\circ T_{\sigma(n)}$. 
For $n \geq m = 3$, define $W_n = \overline{co}(\bigcap_{\sigma \in \mathfrak{S}_n} T^{\sigma}(K))$. 
As $\mathfrak{F}$ is left reversible, $W_n$ is nonempty for all $n \geq 3.$
As for $l = 1, 2, 3$, $T_l(K) \subseteq K$, we have $W_{n+1} \subseteq W_n$ for $n \geq 3$. 

Define $K_0$ be the set of asymptotic center of $\{W_n: n \in  \mathbb{N}\}$ with respect to $K$, where $W_1 = K = W_2$. 


It is claimed that $T_l(K_0) \subseteq K_0$ for $l = 1, 2, 3.$ 
Set $V_{n+1, l} = \overline{co}(\bigcap_{\sigma \in \mathfrak{S}_n}T_l\circ T^{\sigma}(K))$, for $l = 1, 2, 3.$ Note that $W_{n+1} \subseteq V_{n+1, l}$, for $n \geq 3$.

Let $x \in K_0$. Then  
\begin{center}
$\delta(T_l(x), W_{n+1}) \leq \delta(T_l(x), V_{n+1, l})= \delta(x, W_n)$, for $n \geq 3$.
\end{center}
Hence $r(T_l(x)) = \inf \delta(T_l(x), W_n) \leq r(x)$, for all $x \in K_0$. This implies that $T_l(K_0) \subseteq K_0$, for $l = 1, 2, 3$ and consequently  $T(K_0)\subseteq K_0$ for all $T$ in the semigroup $\mathfrak{F}_0$ generated by $\{T_1, T_2, T_3\}$. Therefore, by Theorem \ref{P:T:CFT_2}, the semigroup $\mathfrak{F}_0$ has a common fixed point in $K_0$. Let $x_0 \in K_0$ be a common fixed point of $\mathfrak{F}_0$. 
Now note that for $n \geq 3$
\begin{center}
$\delta(x_0, T^{\sigma}(K)) = \delta(x_0, K)$, for all $\sigma \in \mathfrak{S}_n $. 
\end{center} 
Hence $\delta(x_0, W_n) \leq \delta(x_0, K)$ for all $n \in \mathbb{N}$. 

%
%

Also, note that as every sub-semigruop of $\mathfrak{F}$ is left reversible, there is a $S_n \in \mathfrak{F}_0$ such that $T^{\sigma_0}\circ S_n(K) \subseteq W_n$, where $\sigma_0$ is the identity permutation in $\mathfrak{S}_n$. Since $x_0$ is a common fixed point of $\mathfrak{F}_0$, we have 
\begin{center}
$\|x_0 - T^{\sigma_0}\circ S_N(y)\| = \|x_0 - y\|$, for all $y \in K$.
\end{center}
 Therefore 
\begin{center}
$\delta(x_0, K) = \delta(x_0, T^{\sigma_0}\circ S_N(K)) \leq \delta(x_0, W_N)$, for all $n$.
\end{center}
Hence $r(x_0) = \delta(x_0, K) \leq r(y) \leq \delta(y, K)$, for all $y \in K$, as $x_0 \in K_0$. Consequently, the family $\{T_1, T_2, T_3\}$ has a common fixed in $C(K)$. This proves that the family $\mathfrak{S}$ has the finite intersection property. Therefore, the left reversible semigroup $\mathfrak{F}$ has a common fixed point in $C(K)$.
\hfill{$\Box$}

\begin{remark}
It is known \cite[page-516]{Day_1957} that a subgroup of an amenable group is  amenable. 
Also, note that \cite{Garn_1964,Mitc_1970} every left amenable semigroup is left reversible. Therefore, there are left reversible semigroups in which every sub-semigroup is also left reversible. 
\end{remark}

\begin{theorem}
Let $K$ be a nonempty weakly compact convex set having normal structure in a strictly convex Banach space $X$. Let $\mathfrak{F}$ be a left reversible topological semigroup of isometry, separately continuous actions on $K$. Assume that every sub-semigroup of $\mathfrak{F}$ is also left reversible. Then the Chebyshev center $C(K)$ contains a common fixed point of $\mathfrak{F}$.
\end{theorem}

\noindent{\bf Proof} Note that $F_T=\{x \in C(K): Tx = x\}$ is a non-empty closed convex subset of $C(K)$, whenever $T$ is an isometry on $K$ as $X$ is strictly convex. Hence it is enough to prove that the family $\mathfrak{S} = \{F_T: T \in \mathfrak{F}\}$ has the finite intersection property. 

Consider a finite subset $\{T_1, T_2, ..., T_m\}$ of $\mathfrak{F}$. Take $m = 3$. Define $W_n$ and $K_0$ as in Theorem \ref{P:T:1}. It is easy to see by the arguments in Theorem \ref{P:T:1} that $K_0$ is invariant under every member $T$ in the semigroup $\mathfrak{F}_0$ generated by $\{T_1, T_2, T_3\}$. 

Now, by Theorem \ref{P:T:CFT_2} $K_0$ contains a common fixed point of $\mathfrak{F}_0$. Let $x_0$ be a common fixed point of $\mathfrak{F}_0$. Then  $x_0$ is a common fixed point of $\{T_1, T_2, T_3\}$ and  
\begin{center}
$\delta(x_0, T^{\sigma}(K)) = \delta(x_0, K)$, for all $\sigma$ in the permutation group $\mathfrak{S}_n$. 
\end{center} 
Hence $\delta(x_0, W_n) \leq \delta(x_0, K)$ for all $n \geq 3$. Now, by using the hypothesis that every sub-semigroup of $\mathfrak{F}$ is also left reversible, it can be seen, by using arguments in Theorem \ref{P:T:1}, that $r(x_0) = \inf  \delta(x_0, W_n) = \delta(x_0, K)$. But $r(x_0) \leq r(y) \leq \delta(y, K)$, for all $y \in K$, as $x_0 \in K_0$. This implies that $x_0 \in C(K)$ and consequently the family $\{T_1, T_2, T_3\}$ has a common fixed in $C(K)$. Therefore the family $\mathfrak{S}$ has the finite intersection property.  
This shows that the left reversible semigroup $\mathfrak{F}$ has a common fixed point in $C(K)$.
\hfill{$\Box$}


\end{document}